\documentclass[12pt]{article}
\usepackage{amssymb}

\usepackage{graphicx}
\usepackage{amsmath}

\newtheorem{theorem}{Theorem}

\newtheorem{lemma}[theorem]{Lemma}

\newtheorem{proposition}[theorem]{Proposition}

\begin{document}

\title{Large time behavior of the solutions to the difference wave operators.}
\author{ H. Islami,  B. Vainberg
\thanks {Corresponding author, email: brvainbe@uncc.edu. The work by B. Vainberg was
supported partially by the NSF grant DMS-0405927.}
 \and \it {Dept. of
Mathematics and Statistics, University of NC at Charlotte,} \\
\it {Charlotte, NC 28223.} }
\date{}

\maketitle

\begin{abstract}
The Cauchy problem for two dimensional difference wave operators is
considered with potentials and initial data supported in a bounded region.
The large time asymptotic behavior of solutions is obtained. In contrast to
the continuous case (when the problem in the Euclidian space is considered,
not on the lattice) the resolvent of the corresponding stationary problem
has singularities on the continuous spectrum, and they contribute to the
asymptotics.


\end{abstract}

\textbf{I. Introduction. \ }In the recent years there has been considerable
interest in the scattering theory of discrete Schr\"{o}dinger operator on
the lattice $Z^{d}$( see, for example, \cite{B}-\cite{JMP},\cite{mv1}-\cite
{s}). The two dimensional Schr\"{o}dinger operator on the lattice $Z^{2}$%
\begin{equation*}
Hu=-\Delta u+q\left( \xi \right) u,\text{ \ \ }\xi \in Z^{2}
\end{equation*}
is considered in this paper. Here $\Delta $ is the lattice Laplacian
\begin{equation*}
\Delta u\left( \xi \right) =\sum\limits_{\left| \xi -\xi ^{\prime }\right|
=1}u\left( \xi ^{\prime }\right) -4u\left( \xi \right) .
\end{equation*}
and the potential $q$ is real valued. Our goal is to study the large time $%
\left( t\rightarrow \infty \right) $ asymptotic behavior of solutions of the
Cauchy problem for the wave equation:
\begin{eqnarray}
v_{tt}\left( t,\xi \right) &=&\Delta v\left( t,\xi \right) -q\left( \xi
\right) v\left( t,\xi \right) ,\text{ \ \ \ }t>0,\text{ \ }\xi \in Z^{2};
\notag \\
v|_{t=0} &=&0,\text{ \ \ \ }v_{t}|_{t=0}=f\left( \xi \right) ,  \label{1}
\end{eqnarray}
where the potential $q$ and the initial perturbation $f$ have bounded
supports. Without loss of the generality, one can assume that the function $%
f $ is also real valued.

This question is well studied (see \cite{lp},\cite{v}) in the continuous
case, i.e., for operators $H$ on $L^{2}\left( R^{n}\right) .$ In this case
the answer depends on analytical properties of the truncated resolvent $\hat{%
R}_{\lambda }=\chi \left( H-\lambda \right) ^{-1}\chi $ where $\chi \in
C_{0}^{\infty }.$ In the 2-D continuous case, the truncated resolvent is
analytic in $\lambda $ with the branch point at $\lambda =0$ of the
logarithmic type and with poles at the eigenvalues of $H$. The large time
asymptotic behavior of $v$ in the continuous case is expressed through the
eigenvalues of $H$ and the asymptotics of $\hat{R}_{k^{2}}$ as $k=\sqrt{%
\lambda }\rightarrow 0.$ In fact, we neglect analytic in $k$ terms of the
asymptotic expansion of $\hat{R}_{k^{2}}$ at $k=0.$ So, let $\hat{R}%
_{k^{2}}^{\prime }$ be the operator function $\hat{R}_{k^{2}}$ by modulus of
operator functions in $L^{2}(R^{2})$ which are analytic in $k$ in a
neighborhood of $k=0.$ Then, in the continuous case,
\begin{equation}
\hat{R}_{k^{2}}^{\prime }=Qk^{\alpha }\left( \text{log}k\right) ^{\beta
}+O\left( k^{^{\alpha }}\left( \text{log}k\right) ^{\beta -1}\right) \text{
as }\lambda =k^{2}\rightarrow 0,  \label{2}
\end{equation}
where $Q$ is a bounded (finite dimensional) operator in $L^{2}(R^{2}),$ and $%
\alpha ,\beta $ are integers such that $\alpha \geq -1$ and $\beta \leq -2$
if $\alpha =-1.$

We show that the truncated resolvent $\hat{R}_{\lambda }$ of the two
dimensional lattice Schr\"{o}dinger operator has three logarithmic branch
points at $\lambda =0,\lambda =4$, and $\lambda =8$, and the large time
asymptotics of $v$ can be expressed through eigenvalues of $H$ and the
asymptotics of $\hat{R}_{k^{2}}$ at the branch points. In order to compare
the results valid in the continuous case and in the lattice case considered
below, let us first assume that the operator $H$ does not have eigenvalues.
Then the solution $v\left( t,x\right) $ of the Cauchy problem in the
continuous case has the following asymptotic behavior as $t\rightarrow
\infty $%
\begin{equation}
\chi \left( x\right) v\left( t,x\right) =u\left( x\right) t^{-\alpha
-1}\left( \text{log}t\right) ^{\beta -\gamma }+w\left( t,x\right) ,
\label{3}
\end{equation}
where $u\left( x\right) =cQf,$ $c$ is a constant, $\gamma =0$ if $\alpha <0,$
$\gamma =1$ if $\alpha \geq 0,$ and $\left\| w\right\| \leq Ct^{-\alpha
-1}\left( \text{log}t\right) ^{\beta -\gamma -1}$ as $t\rightarrow \infty .$
In fact if $q\geq 0$ and $q>0$ at least at one point, then $\hat{R}%
_{k^{2}}^{\prime }$ is bounded at $k=0$. Thus $v$ decays at least as $%
t^{-1}\left( \text{log}t\right) ^{-2}$ in this case.

The truncated resolvent $\hat{R}_{k^{2}}^{\prime }$ in the lattice case has
expansions similar to (\ref{2}), when $k=\sqrt{\lambda }\rightarrow 0,$ $%
k\rightarrow \pm 2$ and $k\rightarrow \pm \sqrt{8}$. Correspondingly, if $H$
does not have eigenvalues, then
\begin{eqnarray}
v\left( t,\xi \right) &=&u_{0}\left( \xi \right) t^{-\alpha _{0}-1}\left(
\text{log}t\right) ^{\beta _{0}-\gamma _{0}}+u_{1}\left( \xi \right)
t^{-\alpha _{1}-1}\left( \text{log}t\right) ^{\beta _{1}-\gamma _{1}}\cos
(2t+\omega _{1})  \notag \\
&&u_{2}\left( \xi \right) t^{-\alpha _{2}-1}\left( \text{log}t\right)
^{\beta _{2}-\gamma _{2}}\cos (\sqrt{8}t+\omega _{2})+w(t,\xi ),
\label{asym}
\end{eqnarray}
where $\omega _{j}$ are constants and $w$ at any fixed $\xi $ decays, as $%
t\rightarrow \infty ,$ at least as the biggest term in the right hand side
above multiplied by $\left( \text{log}t\right) ^{-1}$. Note that we can not
guarantee the boundedness of $\hat{R}_{k^{2}}^{\prime }$ in neighborhoods of
$k=\pm 2,$ $k=\pm \sqrt{8}$ in the lattice case when $q\geq 0,$ $q\neq 0$.
Thus, $v$ can decay very slowly even if the potential is nonnegative.

In the continuous case, the operator $H$ may have negative eigenvalues $%
\left\{ \lambda _{j}=-k_{j}^{2}\right\} _{j=1}^{N},$ and then the
exponentially growing terms $w_{j}\left( x\right) e^{k_{j}t}$ have to be
added to the expansions (\ref{3}). This is also true in the lattice case.
The operator $H$ in the continuous case does not have positive eigenvalues
(they can not be embedded into the continuous spectrum). In the lattice
case, the continuous spectrum is the segment $[0,8],$ eigenvalues $\lambda
_{j}\notin \left( 0,4\right) \cup \left( 4,8\right) $, but we cannot exclude
the situation when $\lambda _{j}=4$ or $8.$ Besides, positive eigenvalues $%
\lambda _{j}$ may belong to the ray $(8,\infty ).$ All these positive
eigenvalues contribute additional, oscillating in $t,$ terms $w_{j}\left(
\xi \right) \sin (\sqrt{\lambda _{j}}t)$ into the expansion (\ref{asym}). In
both continuous and lattice cases, the existence of the zero eigenvalue
provides an additional, linear in $t,$ term in (\ref{asym}).

The following approach will be used to prove all the statements
above. From standard a priori estimates for the solutions of
(\ref{1}) it follows that the solutions grow not faster than some
exponent as $t\rightarrow \infty $. Hence, there exists a constant
$A<\infty $ such that the Laplace transform with respect to
$t$-variable $v(t,x)\rightarrow \tilde{v}\left( \mu ,\xi \right)
=\int\limits_{0}^{\infty }v\left( t,\xi \right) e^{-\mu t}dt$
exists for Re$\mu >A$ and satisfies the equation
\begin{equation}
\mu ^{2}\tilde{v}=\Delta \tilde{v}-q\left( \xi \right) \tilde{v}+f\left( \xi
\right) ,\text{ \ \ \ Re}\mu >A,\text{ \ \ \ \ }f,\tilde{v}\in l^{2}(Z^{2}),
\label{6}
\end{equation}
and the inverse transform is given by the integral
\begin{equation*}
v\left( t,\xi \right) =\frac{i}{2\pi }\int_{B-i\infty }^{B+i\infty }\tilde{v}%
\left( \mu ,\xi \right) e^{\mu t}d\mu ,\text{ \ \ }B>A,\text{ \ \ \ }v\left(
t,\cdot \right) \in l^{2}(Z^{2}).
\end{equation*}
Set $\mu =-ik$ and $\tilde{v}\left( -ik,\xi \right) =u\left( k,\xi \right) $%
. Then
\begin{equation}
\left( -\Delta +q-k^{2}\right) u\left( k,\xi \right) =f\left( \xi \right) ,%
\text{ \ \ Im}k>A,\text{ \ \ \ \ }f,u\in l^{2}(Z^{2}),  \label{7}
\end{equation}
\ and
\begin{equation}
v\left( t,\xi \right) =\frac{1}{2\pi }\int_{B-i\infty }^{B+i\infty }u\left(
k,\xi \right) e^{-ikt}dk=\frac{1}{2\pi }\int_{B-i\infty }^{B+i\infty
}R_{k^{2}}fe^{-ikt}dk,\text{ \ \ }B>A,  \label{9}
\end{equation}
where $R_{k^{2}}=(-\Delta +q-k^{2})^{-1}.$ The formula (\ref{9}) is the
staring point for the investigation of the large time asymptotic behavior of
$v$. In order to get the asymptotics we are going to move down the contour
of integration in (\ref{9}), and for this reason we need to know the
analytic properties of the resolvent $R_{k^{2}}.$ The next section is
devoted to a study of the analytic properties of the resolvent $%
R_{k^{2}}^{0} $ of the unperturbed operator $H=-\Delta .$ The properties of
the $R_{k^{2}}$ will be studied in the section 3, and the last section deals
with the large time asymptotics of $v$.

\textbf{II. Unperturbed Problem. }If $q\left( \xi \right) \equiv 0$ then (%
\ref{7}) becomes
\begin{equation}
\left( -\Delta -k^{2}\right) u\left( k,\xi \right) =f\left( \xi \right) .
\label{10}
\end{equation}

We solve (\ref{10}) using the Fourier transform with respect to $\xi $%
-variable. We denote by $T$ the open square $(-\pi ,\pi )\times (-\pi ,\pi )$
in $R^{2}$. Since $f$ has bounded support, there is a square $S\subset Z^{2}$
which contains the support of $f$. Then the Fourier transform of $f$ is a
function on $T:$%
\begin{equation}
\hat{f}\left( \sigma \right) =\sum\limits_{\xi \in S}f\left( \xi \right)
e^{-i\sigma \xi },\text{ \ \ }\sigma =\left( \sigma _{1},\sigma _{2}\right)
\in T,  \label{11}
\end{equation}
where $\sigma \xi =\sigma _{1}\xi _{1}+\sigma _{2}\xi _{2}.$

The Fourier transform of the lattice operator $-\Delta $ is the operator of
multiplication by
\begin{equation*}
\varphi \bigskip \left( \sigma \right) =4-2\cos \sigma _{1}-2\cos \sigma
_{2},
\end{equation*}
and the function $\varphi $ is real valued and its range is $\left[ 0,8%
\right] $. Thus, $-\Delta $ is a self-adjoint operator with the spectrum $%
\left[ 0,8\right] ,$ and the resolvent $R_{\lambda }=\left( -\Delta -\lambda
\right) ^{-1}$ is analytic in the complex $\lambda $-plane with the cut
along the interval $\left[ 0,8\right] $. It is convenient for us to use the
variable $k=\sqrt{\lambda }$. If $k$ is considered as a spectral parameter,
then the upper half plane $C_{+}=\{k:$ Im$k>0\}$ belongs to the resolvent
set, and the spectrum is $\left[ -\sqrt{8},\sqrt{8}\right] $. Hence, from (%
\ref{7}) it follows that
\begin{equation}
u\left( k,\xi \right) =R_{k^{2}}f=\frac{1}{2\pi }\int\limits_{T}\frac{\hat{f}%
\left( \sigma \right) e^{i\sigma \xi }}{\varphi \left( \sigma \right) -k^{2}}%
d\sigma ,\text{ \ \ \ \ }k\in \overline{C}_{+}\backslash \left[ -\sqrt{8},%
\sqrt{8}\right] .  \label{13}
\end{equation}
Formulas (\ref{13}) and (\ref{11}) imply that
\begin{equation}
u\left( k,\xi \right) =\frac{1}{2\pi }\sum\limits_{\eta \in S}G\left( k,\xi
-\eta \right) f\left( \eta \right) ,\text{ \ \ \ \ }k\in \overline{C}%
_{+}\backslash \left[ -\sqrt{8},\sqrt{8}\right] .  \label{14}
\end{equation}
where
\begin{equation}
G=G\left( k,\xi \right) =\frac{1}{2\pi }\int\limits_{T}\dfrac{e^{i\sigma \xi
}}{\varphi \left( \sigma \right) -k^{2}}d\sigma ,\ \ \ \ k\in \overline{C}%
_{+}\backslash \left[ -\sqrt{8},\sqrt{8}\right] ,  \label{gre}
\end{equation}
is Green's function of the operator $-\Delta -k^{2}.$ Let us denote by $%
k_{s} $ the following points: $k_{0}=0,k_{\pm 1}=\pm 2$ and $k_{\pm 2}=\pm
\sqrt{8} $.

\begin{theorem}
\label{t1}. For each fix $\xi \in Z^{2},$ Green's function $G$ is analytic
in $k$ when $k\in \overline{C}_{+}\backslash \left[ -\sqrt{8},\sqrt{8}\right]
$ and it admits an analytic extension on $\left[ -\sqrt{8},\sqrt{8}\right]
\backslash \cup k_{s}.$ Green's function $G$ in a neighborhood of each point
$k_{s}$ has the form
\begin{equation}
G\left( k,\xi \right) =u_{1,s}(k,\xi )\text{log}\left( k-k_{s}\right)
+u_{2.s}\left( k,\xi \right)  \label{16}
\end{equation}
where $u_{1,s},u_{2,s}$ are analytic in $k$ for each $\xi \in Z^{2}.$\bigskip
\end{theorem}

The proof of this Theorem will be based on the following three Lemmas. Since
we are going to study the function $G$ when $\xi $ is fixed, we will often
omit the variable $\xi $ from the argument of $G$ (and some other functions)
and write simply $G=G(k).$

Let $D$ be a bounded domain in $R_{\sigma }^{2},$ $\sigma =(\sigma
_{1},\sigma _{2}),$ with a piece-wise analytic boundary. Let $\psi =\psi
(\sigma )$ be a real valued, analytic in $\sigma \in \overline{D}$ function
and
\begin{equation*}
\Gamma =\{\sigma \in \overline{D}:\text{ }\psi =0\}.
\end{equation*}
Consider
\begin{equation}
h(z)=\int_{D}\frac{r(\sigma )}{\psi (\sigma )-z}d\sigma ,\text{ \ \ \ Im}z>0,
\label{ah}
\end{equation}
where the function $r$ is continuous in $\overline{D}$ and analytic in a
neighborhood of $\Gamma .$ Obviously, $h(z)$ is analytic in the half plane Im%
$z>0.$

\begin{lemma}
\label{ll1}. Let $\nabla \psi \neq 0$ on $\Gamma $, and\ $\Gamma \cap
\partial D$ be empty or belong to an analytic part of $\partial D$ with the
intersection being transversal. Then $h(z)$ admits an analytic extension
into a neighborhood of the point $z=0.$
\end{lemma}

\textbf{Proof}. One may assume that $\Gamma $ is connected, since otherwise $%
D$ can be split into parts containing only connected components of $\Gamma .$
Let $\Gamma ^{\prime }$ be a small analytic extension of $\Gamma $ beyond $%
\overline{D}$ if $\Gamma \cap \partial D$ is not empty, and let $\Gamma
^{\prime }=\Gamma $ if $\Gamma \cap \partial D$ is empty. Let $\sigma
=\sigma (s),$ $s\in \lbrack 0,L],$ be an analytic parametrization of $\Gamma
^{\prime }.$ Consider the following problem
\begin{equation}
\frac{d\sigma }{dt}=\frac{\nabla \psi (\sigma )}{|\nabla \psi (\sigma )|^{2}}%
,\text{ \ \ \ \ }\sigma |_{t=0}=\sigma (s),\text{ \ \ }s\in \lbrack 0,L].
\label{par}
\end{equation}
Solution $\sigma =\sigma (t,s)$ of (\ref{par}) exists and is an analytic
function of $t$ and $s$ when $|t|$ is small enough. Since the Jacobian $J=%
\frac{d\sigma }{dtds}$ is not zero on $\Gamma ^{\prime },$ there is an $%
\varepsilon >0$ such that $\sigma =\sigma (t,s)$ is analytic and $J\neq 0$
on
\begin{equation}
\Gamma _{\varepsilon }^{\prime }=\{\sigma =\sigma (t,s):\text{ }s\in \lbrack
0,L],\text{ }|t|\leq \varepsilon \}.  \label{gap}
\end{equation}

We split $D$ into two non-intersecting parts: $D=D_{1}\cup D_{2},$ where $%
D_{1}=D\cap \Gamma _{\varepsilon }^{\prime },$ and we represent $h(z)$ as a
sum
\begin{equation}
h(z)=h_{1}(z)+h_{2}(z),\text{ \ \ \ }h_{i}(z)=\int_{D_{i}}\frac{r(\sigma )}{%
\psi (\sigma )-z}d\sigma .  \label{hi}
\end{equation}

Since $\psi (\sigma )\neq 0$ on $\overline{D}_{2},$ we have $|\psi |>\delta
>0$ on $\overline{D}_{2},$ and therefore $h_{2}$ is analytic when $%
|z|<\delta .$ In order to study $h_{1}$ we change the variables $\sigma
\rightarrow (t,s).$ Note that
\begin{equation*}
\frac{d\psi (\sigma )}{dt}=\nabla \psi (\sigma )\cdot \frac{d\sigma }{dt}=1,%
\text{ \ \ }\psi (\sigma )|_{t=0}=0.
\end{equation*}
Thus, $\psi (\sigma )=t,$ and formula (\ref{hi}) for $h_{1}$ takes the form
\begin{equation}
h_{1}(z)=\int_{D_{1}^{\prime }}\frac{v(t,s)}{t-z}dsdt,  \label{h1}
\end{equation}
where $D_{1}^{\prime }$ is the image of $D$ under the map $\sigma
\rightarrow (t,s)$ and $v=r(\sigma )J$ is analytic in $(t,s)\in
D_{1}^{\prime }$ if $\varepsilon $ is small enough. Since $\Gamma
_{\varepsilon }^{\prime }$ is defined by inequalities $0\leq s\leq L,$ $%
|t|\leq \varepsilon $ (see(\ref{gap})), (\ref{h1}) when $\varepsilon $ is
small enough can be rewritten in the form
\begin{equation*}
h_{1}(z)=\int_{-\varepsilon }^{\varepsilon }\int_{s_{1}(t)}^{s_{2}(t)}\frac{%
v(t,s)}{t-z}dsdt,
\end{equation*}
where the functions $s_{i}(t)$ are analytic. After the integration with
respect to the variable $s$ we get
\begin{equation*}
h_{1}(z)=\int_{-\varepsilon }^{\varepsilon }\frac{w(t)}{t-z}dt,
\end{equation*}
where $w$ is analytic. Thus
\begin{equation*}
h_{1}(z)=\int_{-\varepsilon }^{\varepsilon }\frac{w(t)-w(z)}{t-z}%
dt+w(z)\int_{-\varepsilon }^{\varepsilon }\frac{1}{t-z}dt.
\end{equation*}
The first integrand above (defined as $w^{\prime }(z)$ when $t=z$) is
analytic in $t$ and $z$ when $|t|$ and $|z|$ are small$,$ and therefore the
first integral is analytic in a neighborhood of $z=0.$ The second integral
is equal to $w(z)$log$\frac{\varepsilon -z}{-\varepsilon -z}.$ So, it is
also analytic. Lemma \ref{ll1} is proved.

\begin{lemma}
\label{l0}. Let function $h$ be given by (\ref{ah}) with $\Gamma $
consisting of one point $\sigma =\sigma ^{0}\in D,$ and let $\det A>0$ where
$A=[\frac{\partial ^{2}\psi }{\partial \sigma _{i}\partial \sigma _{j}}%
(\sigma ^{0})].$ Then $h(z)$ has the following form in a neighborhood of the
point $z=0$%
\begin{equation}
h(z)=u_{1}(z)\text{log}z+u_{2}(z),  \label{hlna}
\end{equation}
where $u_{i}(z)$ are analytic in a neighborhood of the origin, and
\begin{equation}
u_{1}(0)=-\pi \gamma r(\sigma ^{0})|\det A|^{-1/2},  \label{111}
\end{equation}
where $\gamma =1$ if eigenvalues of $A$ are positive and $\gamma =-1$ if
they are negative.
\end{lemma}

\textbf{Proof.} We shall prove the Lemma assuming that $A$ has positive
eigenvalues. The other case can be easily reduced to the one which we are
going to consider. From the assumptions of Lemma \ref{ll2} it follows that
in a neighborhood of $\sigma ^{0}$ there exist local coordinates $\alpha
=(\alpha _{1},\alpha _{2}),$ $|\alpha |<\varepsilon ,\ $such that the
vector-function $\sigma =\sigma (\alpha ),$ $|\alpha |\leq \varepsilon ,$ is
analytic, the Jacobian $J=\frac{d\sigma }{d\alpha }\neq 0$ when $|\alpha
|\leq \varepsilon ,$ $J=|\det A|^{-1/2}$ at the point $\sigma ^{0},$ and the
function $\psi $ in the new coordinates is equal to $|\alpha |^{2}.$ \ Let
\begin{equation*}
D_{1}=\{\sigma :\sigma =\sigma (\alpha ),\text{ \ }|\alpha |<\varepsilon \}
\end{equation*}
We take $\varepsilon $ so small that $\overline{D}_{1}\subset D$ and $%
u(\sigma )$ is analytic in $\overline{D}_{1}.$ \ We represent $h(z)$ as a
sum $h_{1}(z)+h_{2}(z)$ where $h_{1},$ $h_{2}$ are the integrals (\ref{ah})
over $D_{1}$ and $D\backslash D_{1},$ respectively. Since $\varphi (\sigma
)-z\neq 0$ when $\sigma \in \overline{D\backslash D}_{1}$ and $|z|$ is small
enough, the function $h_{2}$ is analytic in a neighborhood of $z=0$. Thus it
remains to prove (\ref{hlna}) for the function $h_{1}$.

The change of the coordinates $\sigma \rightarrow \alpha $ allows us to
rewrite $h_{1}$ in the form
\begin{equation}
h_{1}(z)=\int_{|\alpha |<\varepsilon }\frac{v(\alpha )}{|\alpha |^{2}-z}%
d\alpha ,\text{ \ \ \ Im}z>0,  \label{h11}
\end{equation}
where the function $v$ is analytic and
\begin{equation}
v(0)=r(\sigma ^{0})J(\sigma ^{0})=r(\sigma ^{0})|\det A|^{-1/2}.  \label{00}
\end{equation}
Let $v_{1}$ and $v_{2}$ be even and odd parts of $v$:
\begin{equation}
v_{1}(\sigma )=\frac{v(\sigma )+v(-\sigma )}{2},\text{ \ \ \ }v_{2}(\sigma )=%
\frac{v(\sigma )-v(-\sigma )}{2}.  \label{nn}
\end{equation}
We substitute $v_{1}+v_{2}$ for $v$ in (\ref{nn}). Then $%
h_{1}=h_{1,1}+h_{1,2}$ where the terms in the right hand side are given by (%
\ref{h11}) with $v_{1}$ and $v_{2},$ respectively, instead of $v.$
Obviously, $h_{1,2}=0,$ i.e. $h_{1}=h_{1,1}.$ Let $(\rho ,\theta )$ be the
polar coordinates in the $\alpha $ plane. We write the integral for $h_{1,1}$
in the polar coordinates and integrate with respect to $\theta .$ Since $%
v_{1}$ depends analytically on $\rho ^{2},$ we arrive at
\begin{equation}
h_{1}(z)=\int_{0}^{\varepsilon }\frac{\rho w(\rho ^{2})}{\rho ^{2}-z}d\rho ,%
\text{ \ \ \ Im}z>0,  \label{224}
\end{equation}
where the function $w$ is analytic and
\begin{equation}
w(0)=2\pi v(0)=2\pi r(\sigma ^{0})|\det A|^{-1/2}.  \label{v0}
\end{equation}
Formula (\ref{224}) implies that
\begin{equation*}
h_{1}(z)=\int_{0}^{\varepsilon }\frac{\rho \lbrack w(\rho ^{2})-w(z)]}{\rho
^{2}-z}d\rho +w(z)\int_{0}^{\varepsilon }\frac{\rho }{\rho ^{2}-z}d\rho .
\end{equation*}
The first integrand above is analytic in $\rho $ and $z$ when $\rho $ and $%
|z|$ are small$,$ and therefore the first integral is an analytic function
in a neighborhood of $z=0.$ The second integral is equal to $\frac{1}{2}w(z)$%
log$(\varepsilon ^{2}-z)-\frac{1}{2}w(z)$log$(-z).$ The first term in the
last sum is analytic when $|z|$ is small enough, and the second term is a
sum of an analytic function $-\frac{\pi i}{2}w(z)$ and the function $-\frac{1%
}{2}w(z)$log$z.$ This justifies (\ref{hln}) with $u_{1}(z)=-\frac{1}{2}w(z).$
The proof of the Lemma \ref{ll2} is completed

\begin{lemma}
\label{ll2}. Let all the requirements of Lemma \ref{ll1} be satisfied, but $%
\nabla \psi =0$ at one point $\sigma ^{0}$ of $\Gamma $. Let $\sigma ^{0}$ \
belong to the interior part of $D$, and $\det A<0$ where $A=[\frac{\partial
^{2}\psi }{\partial \sigma _{i}\partial \sigma _{j}}(\sigma ^{0})]$. Then
function $h(z),$ defined by (\ref{ah}), has the following form in a
neighborhood of the point $z=0$%
\begin{equation}
h(z)=u_{1}(z)\text{log}z+u_{2}(z),  \label{hln}
\end{equation}
where $u_{i}(z)$ are analytic in a neighborhood of the origin, and $%
u_{1}(0)=i\pi u(\sigma ^{0})|\det A|^{-1/2}.$
\end{lemma}

\textbf{Proof.} There exist analytic local coordinates $\alpha =(\alpha
_{1},\alpha _{2})$ in a neighborhood of $\sigma ^{0}$ in which the function $%
\psi $ has the form $\psi =$ $\alpha _{1}^{2}-\alpha _{2}^{2}.$ Let
\begin{equation*}
D_{1}=\{\sigma =\sigma (\alpha ):|\alpha _{1}|<\varepsilon ,\text{ }|\alpha
_{2}|<2\varepsilon \}
\end{equation*}

We represent $h(z)$ as a sum of two terms $h_{1}(z)+h_{2}(z)$ where $h_{2}$
is the integral over $D\backslash D_{1}$ which is analytic in a neighborhood
of $z=0,$ and the function $h_{1}$ is the integral over $D_{1}$ which has
the following form
\begin{equation}
h_{1}(z)=\int_{-\varepsilon }^{\varepsilon }\int_{-2\varepsilon
}^{2\varepsilon }\frac{v(\alpha )}{\alpha _{1}^{2}-\alpha _{2}^{2}-z}d\alpha
_{2}d\alpha _{1},~~\rm{Im}z>0,  \label{a1}
\end{equation}
where the function $v$ is analytic and (\ref{00}) holds. Lemma \ref{ll1}
implies that the function $h_{2}$ is analytic in a neighborhood of $z=0.$
Hence, it remains to prove the validity of Lemma \ref{ll2} for the function (%
\ref{a1}). We represent the function $v$ as the sum $v_{1}+v_{2}$ where the
function $v_{1}$ is even with respect to $\alpha _{1}$and $v_{2}$ is odd.
Obviously, the integral with the function $v_{2}$ instead of $v$ is equal to
zero. Then we represent $v_{1}$ as the sum of an even with respect to $%
\alpha _{2}$ function and an odd with respect to $\alpha _{2}$ function.
Similarly, only the even part remains. Thus, (\ref{a1}) can be rewritten in
the form
\begin{equation}
h_{1}(z)=\int_{-\varepsilon }^{\varepsilon }\int_{-2\varepsilon
}^{2\varepsilon }\frac{w(\alpha _{1}^{2},\alpha _{2}^{2})}{\alpha
_{1}^{2}-\alpha _{2}^{2}-z}d\alpha _{1}d\alpha _{2}=4\int_{0}^{\varepsilon
}\int_{0}^{2\varepsilon }\frac{w(\alpha _{1}^{2},\alpha _{2}^{2})}{\alpha
_{1}^{2}-\alpha _{2}^{2}-z}d\alpha _{1}d\alpha _{2},\text{ \ \ \ Im}%
z>0,  \label{a2}
\end{equation}
where $w$ is analytic and
\begin{equation*}
w(0)=r(\sigma ^{0})|\det A|^{-1/2}.
\end{equation*}

From (\ref{a2}) it follows that
\begin{equation*}
h_{1}(z)=4\int_{0}^{\varepsilon }\int_{0}^{2\varepsilon }\frac{w(\alpha
_{1}^{2},\alpha _{2}^{2})-w(\alpha _{2}^{2}+z,\alpha _{2}^{2})}{\alpha
_{1}^{2}-\alpha _{2}^{2}-z}d\alpha _{1}d\alpha _{2}+4\int_{0}^{\varepsilon
}\int_{0}^{2\varepsilon }\frac{w(\alpha _{2}^{2}+z,\alpha _{2}^{2})}{\alpha
_{1}^{2}-\alpha _{2}^{2}-z}d\alpha _{1}d\alpha _{2}.
\end{equation*}
The first integrand above is analytic in $\alpha $ and $z,$ and therefore
the first integral is an analytic function in a neighborhood of $z=0.$ \
Thus the statement of the Lemma has to be proved for the second term $g$
which, after integrating over $\alpha _{1}$ and replacing $\alpha _{2}$ by $%
\tau $, takes the form
\begin{equation}
g=2\int_{0}^{\varepsilon }\frac{w(\tau ^{2}+z,\tau ^{2})}{\sqrt{\tau ^{2}+z}}%
\text{log}\frac{\alpha _{1}-\sqrt{\tau ^{2}+z}}{\alpha _{1}+\sqrt{\tau ^{2}+z%
}}|_{\alpha _{1}=0}^{\alpha _{1}=2\varepsilon }d\tau ,~~\rm{Im}z>0. \label{23}
\end{equation}
If $\varkappa =\dfrac{\alpha _{1}-\sqrt{\tau ^{2}+}z}{\alpha _{1}+\sqrt{\tau
^{2}+z}}$ and Im$z>0,$ $\alpha _{1}>0,$ then Im$\varkappa <0.$ The unique
branch of the log$\varkappa $ in (\ref{23}) is chosen by the condition that $%
\pi <\arg ($log$\varkappa )<2\pi \ $\ when Im$z>0$. Thus when Im$z>0,$ we
have
\begin{equation}
g=2\int_{0}^{\varepsilon }\frac{w(\tau ^{2}+z,\tau ^{2})}{\sqrt{\tau ^{2}+z}}%
\text{log}\frac{2\varepsilon -\sqrt{\tau ^{2}+z}}{2\varepsilon +\sqrt{\tau
^{2}+z}}d\tau -2\pi i\int_{0}^{\varepsilon }\frac{w(\tau ^{2}+z,\tau ^{2})}{%
\sqrt{\tau ^{2}+z}}d\tau .  \label{24}
\end{equation}
One can easily check that the function $\dfrac{1}{u}$log$\dfrac{2\varepsilon
-u}{2\varepsilon +u}$ is analytic in $u$ in the circle $\left| u\right|
<2\varepsilon $ (if the function is defined as zero at $u=0$). Since the
integrand in the first term in the right-hand side of (\ref{24}) is analytic
in $\tau $ and $z$ when $\left| z\right| <\varepsilon ^{2}$, we obtain that
the first term in the right-hand side of (\ref{24}) is analytic in $z,$ $%
\left| z\right| <\varepsilon ^{2}$. Hence, it remains to prove the statement
of the Lemma for the function
\begin{equation}
g_{1}\left( z\right) =-2\pi i\int_{0}^{\varepsilon }\frac{w(\tau ^{2}+z,\tau
^{2})}{\sqrt{\tau ^{2}+z}}d\tau ,\text{ \ \ Im}z>0.  \label{25a}
\end{equation}

Obviously $g_{1}$ is an analytic function of $z$ when $z\in U_{\varepsilon
}=\{z:|z|<\varepsilon ^{2},$ $z\notin \overline{R}_{-}\}$, where $R_{-}$ is
the negative part of the real axis. Let us compare the limit values of $%
g_{1} $ when $z\rightarrow -t\pm i0,$ $\ \ 0<\sqrt{t}<\varepsilon .$ We
represent $g_{1}$ in the form $g_{1}\left( z\right) =f_{1}\left( z\right)
+f_{2}\left( z\right) $, where
\begin{equation*}
f_{1}\left( z\right) =-2\pi i\int_{0}^{\sqrt{t}}\frac{w(\tau ^{2}+z,\tau
^{2})}{\sqrt{\tau ^{2}+z}}d\tau ,\text{ \ \ \ }f_{2}\left( \varkappa \right)
=-2\pi i\int_{\sqrt{t}}^{\varepsilon }\frac{w(\tau ^{2}+z,\tau ^{2})}{\sqrt{%
\tau ^{2}+z}}d\tau ,\ \ 0<\sqrt{t}<\varepsilon .
\end{equation*}
Then
\begin{equation*}
f_{1}\left( -t+i0\right) -f_{1}\left( -t-i0\right) =-4\pi \int_{0}^{\sqrt{t}}%
\frac{w(\tau ^{2}-t,\tau ^{2})}{\sqrt{t-\tau ^{2}}}d\tau ,
\end{equation*}
and
\begin{equation*}
f_{2}\left( -t+i0\right) -f_{2}\left( -t-i0\right) =0.
\end{equation*}
Thus, after the substitution $\tau =\sqrt{t}s$, we get
\begin{equation}
g_{1}\left( -t+i0\right) -g_{1}\left( -t+i0\right) =-4\pi \int_{0}^{\sqrt{t}}%
\frac{w(\tau ^{2}-t,\tau ^{2})}{\sqrt{t-\tau ^{2}}}d\tau =-4\pi \int_{0}^{1}%
\frac{w\left( ts^{2}-t,ts^{2}\right) }{\sqrt{1-s^{2}}}ds.  \label{26}
\end{equation}

We denote the latter function by $\varphi _{1}=\varphi _{1}(t).$ Obviously, $%
\varphi _{1}(z)$ is analytic in $z.$ Define
\begin{equation}
\varphi _{2}\left( z\right) =g_{1}\left( z\right) -\frac{1}{2\pi i}\varphi
_{1}(z)\text{log}z.  \label{27}
\end{equation}
Then
\begin{equation*}
\varphi _{2}\left( -t+i0\right) -\varphi _{2}\left( -t-i0\right) =0,
\end{equation*}
i.e. $\varphi _{2}\left( z\right) $ does not have a branch point at $z=0$.
Now in order to prove that $\varphi _{2}$ is analytic at $z=0$ \ it is
enough to show that $\left| \varphi _{2}\left( z\right) \right| \leq c$log$%
\dfrac{1}{\left| z\right| }$ when $\left| z\right| $ is small enough. \ From
(\ref{27}) it follows that it is enough to get the corresponding logarithmic
estimate for $g_{1}.$ We represent $g_{1}$ as the sum of two terms $%
g_{1}=g_{1,1}+g_{1,2}$ by writing the function $w$ in the integral (\ref{25a}%
) as $w_{1}+w_{2}$ where $w_{1}=w(0,-z),$ $w_{2}=w(\tau ^{2}+z,\tau
^{2})-w(0,-z).$ Since
\begin{equation*}
|w_{2}|=|w(\tau ^{2}+z,\tau ^{2})-w(0,-z)|\leq C|\tau ^{2}+z|\text{ \ \ \
for small }|\tau ^{2}+z|,
\end{equation*}
the function $|g_{1,2}|$ is bounded when $|z|$ is small enough. The integral
defining $g_{1,1}$ can be easily evaluated:
\begin{equation*}
g_{1,1}=-2\pi iw(0,-z)\text{log}\frac{\varepsilon +\sqrt{\varepsilon ^{2}+z}%
}{\sqrt{z}}.
\end{equation*}
Hence, $g_{1}$ has the logarithmic estimate as $|z|\rightarrow 0,$ and
therefore $g_{1}\left( z\right) $ is analytic at $z=0.$ The latter, together
with (\ref{27}), proves (\ref{hln}) for $g_{1}$ with $u_{1}(z)=\frac{1}{2\pi
i}\varphi _{1}(z).$ Thus (\ref{ah}) is valid for $h$ with
\begin{equation*}
u_{1}(0)=2iw(0)\int_{0}^{1}\frac{ds}{\sqrt{1-s^{2}}}=\pi iw(0)=\pi ir(\sigma
^{0})|\det A|^{-1/2}.
\end{equation*}
The proof of the Lemma is completed.

\textbf{Proof of the Theorem \ref{t1}}. Let us fix an arbitrary point $%
k^{0}\in \left[ -\sqrt{8},\sqrt{8}\right] \backslash \cup k_{s}$ and put $%
k=k^{0}+\varkappa $ in (\ref{gre}). Then the function $2\pi G$ defined by (%
\ref{gre}) takes the form (\ref{ah}) with $\psi (\sigma )=\varphi (\sigma
)-(k^{0})^{2}$ and $z=2k^{0}\varkappa +\varkappa ^{2}.$ One can easily check
that the integral (\ref{ah}), which we got, satisfies all the assumptions of
Lemma \ref{ll1} (with $\Gamma $ being a closed curve). Thus $%
G(k^{0}+\varkappa ,\xi )$ is analytic in $z$ when $|z|$ is small enough, and
therefore it is analytic in $\varkappa $ in a neighborhood of the point $%
\varkappa =0.$ This proves the first statement of Theorem \ref{t1}.

The statement of Theorem \ref{t1} concerning the point $k=k_{0}=0$ is the
direct consequence of Lemma \ref{l0}, since $2\pi G$ defined by (\ref{gre})
has the form (\ref{ah}) with $\psi (\sigma )=\varphi (\sigma )$ and $%
z=k^{2}, $ and all the assumptions of Lemma \ref{l0} hold in this case (with
$\sigma ^{0}=0)$. In order to prove the statement concerning the points $%
k=k_{\pm 2}, $ we use the periodicity of the integrand in (\ref{gre}) and
replace the square $T$ in (\ref{gre}) by $T^{\prime }=(0,2\pi )\times
(0,2\pi )\subset R_{\sigma }^{2}.$ After that, we again can apply Lemma \ref
{l0} with $\psi (\sigma )=\varphi (\sigma )-8,$ $z=$ $k^{2}-8$ and $\sigma
^{0}=(\pi ,\pi ).$ In order to prove the statement of the Theorem concerning
the points $k=k_{\pm 1},$ we replace the square $T$ in (\ref{gre}) by $%
T^{\prime \prime }=(\frac{-\pi }{2},\frac{3\pi }{2})\times (\frac{-\pi }{2},%
\frac{3\pi }{2})$ and split $T^{\prime \prime }$ in two parts by the line $%
\sigma _{1}=\sigma _{2}.$ Lemma \ref{ll2} can be applied to each of these
two integrals with $\psi (\sigma )=\varphi (\sigma )-4,$ $z=$ $k^{2}-4$ and
the point $\sigma ^{0}$ equal to $(0,\pi )$ in the first integral, and equal
to $(\pi ,0)$ in the second integral.

The proof of the Theorem is completed.

Later, when the large time asymptotic behavior of the solutions to the
Cauchy problem is analyzed, we shall need to know some properties of the
coefficients $u_{1,s}$ and $u_{2,s}$ in (\ref{16}). To be more exact, we
will need the following

\begin{lemma}
\label{l00}. Functions $u_{1,0}\left( 0,\xi \right) $ and $u_{2,0}\left(
0,\xi \right) \ $have the following properties:

a) $u_{1,0}\left( 0,\xi \right) =c_{1}$ is a constant,

b) $u_{2,0}\left( 0,\xi \right) =c_{1}\ln \left| \xi \right| +c_{2}+o\left(
\left| \xi \right| ^{-1}\right) $ as $\left| \xi \right| \rightarrow \infty
. $
\end{lemma}

\textbf{Proof.} Let $\varsigma \in C_{0}^{\infty }(R_{\sigma }^{2}),$ and $%
\varsigma \left( \sigma \right) =1$ if $\left| \sigma \right| <1,$ $%
\varsigma \left( \sigma \right) =0$ if $\left| \sigma \right| >2.$ From (\ref
{gre}) it follows that

\begin{equation}
2\pi G\left( k,\xi \right) =\int_{R^{2}}\varsigma \left( \sigma \right)
\frac{e^{i\sigma \xi }}{\varphi \left( \sigma \right) -k^{2}}d\sigma
+\int_{T}\left( 1-\varsigma \left( \sigma \right) \right) \frac{e^{i\sigma
\xi }}{\varphi \left( \sigma \right) -k^{2}}d\sigma .  \label{a45}
\end{equation}
We denote the second term in the right hand side of (\ref{a45}) by $%
v=v(k,\xi ).$ Since $\varphi \left( \sigma \right) =0$ only at $\sigma =0,$
and $1-\varsigma =0$ in a neighborhood of the origin, then $v$ is analytic
in $k$ for small enough $\left| k\right| ,$ and
\begin{equation*}
v=v(0,\xi )=\int_{T}\left( 1-\varsigma \left( \sigma \right) \right) \frac{%
e^{i\sigma \xi }}{\varphi \left( \sigma \right) }d\sigma .
\end{equation*}
The integrand above is infinitely smooth. We apply the integration by parts
to that integral, integrating $e^{i\sigma \xi }$ and differentiating the
factor $\dfrac{1-\varsigma \left( \sigma \right) }{\varphi \left( \sigma
\right) }$. Since $1-\varsigma $ is equal to one in a neighborhood of the
boundary of $T$ and the other factors in the integrand are periodic, the
contributions from the boundary will be cancelled, and the integration by
parts leads to the estimate:
\begin{equation}
v\left( 0,\xi \right) =O(\left| \xi \right| ^{-\infty })\text{ as }\left|
\xi \right| \rightarrow \infty .  \label{a46}
\end{equation}

Consider
\begin{equation*}
w\left( k,\xi \right) =\int_{R^{2}}\varsigma \left( \sigma \right) \frac{%
e^{i\sigma \xi }}{\left| \sigma \right| ^{2}-k^{2}}d\sigma .
\end{equation*}

Let us study
\begin{equation}
2\pi G-v-w=\int_{R^{2}}\varsigma e^{i\sigma \xi }\left[ \frac{1}{\varphi
\left( \sigma \right) -k^{2}}-\frac{1}{\left| \sigma \right| ^{2}-k^{2}}%
\right] d\sigma =\int_{R^{2}}\varsigma e^{i\sigma \xi }\frac{\left| \sigma
\right| ^{2}-\varphi \left( \sigma \right) }{\left( \varphi \left( \sigma
\right) -k^{2}\right) (\left| \sigma \right| ^{2}-k^{2})}d\sigma .
\label{a48}
\end{equation}
The numerator of the integrand in (\ref{a48}) is of order $O(\left| \sigma
\right| ^{4})$ as $\left| \sigma \right| \rightarrow 0$. Thus, the integral
has a limit as $k\rightarrow 0:$%
\begin{equation}
\left( 2\pi G-v-w\right) \left( 0,\xi \right) =\int_{R^{2}}e^{i\sigma \xi
}\psi \left( \sigma \right) d\sigma ,\text{ \ \ \ where }\psi \left( \sigma
\right) =\varsigma \frac{\left| \sigma \right| ^{2}-\varphi \left( \sigma
\right) }{\varphi \left( \sigma \right) \left| \sigma \right| ^{2}}.
\label{a49}
\end{equation}
Since $\nabla \psi \left( \sigma \right) \in L^{1},$ its Fourier transform
is $o|\xi |^{-1}$as $|\xi |\rightarrow \infty .$ From here, (\ref{a49}), and
(\ref{a46}) it follows that
\begin{equation}
\left| \left( 2\pi G-w\right) \left( 0,\xi \right) \right| =o(\left| \xi
\right| ^{-1})\text{ as }\left| \xi \right| \rightarrow \infty .  \label{a51}
\end{equation}

Consider
\begin{equation*}
w_{1}\left( k,\xi \right) =\int_{R^{2}}\left( 1-\varsigma \right) \frac{%
e^{i\sigma \xi }}{\left| \sigma \right| ^{2}-k^{2}}d\sigma ,~~\rm{ Im}k>0.
\end{equation*}
Then $w_{1}$ is analytic in $k$ for $\left| k\right| <\varepsilon $ and
\begin{equation}
w_{1}\left( 0,\xi \right) =O\left( \left| \xi \right| ^{-\infty }\right)
\text{ as }\left| \xi \right| \rightarrow \infty .  \label{a52}
\end{equation}
Thus,
\begin{equation}
|\left( 2\pi G-\left( w+w_{1}\right) \right) |\left( 0,\xi \right) =O(\left|
\xi \right| ^{-1})\text{ as }\left| \xi \right| \rightarrow \infty .
\label{a53}
\end{equation}
Note that $(w+w_{1})/2\pi $ is the inverse Fourier transform of the function
$\dfrac{1}{\left| \sigma \right| ^{2}-k^{2}}$ in $R_{\sigma }^{2}$, i.e.
\begin{equation*}
w+w_{1}=c_{0}H_{0}\left( k\left| \xi \right| \right) ,
\end{equation*}
where $c_{0}$ is a constant, and $H_{0}$ is the Hankel function. Hence,

\begin{equation}
w+w_{1}=2\pi \lbrack c_{1}\text{log}\left( k\left| \xi \right| \right)
+c_{2}]+\alpha \left( k\left| \xi \right| \right) ,  \label{a55}
\end{equation}
where the function $\alpha $ is continuous and $\alpha (0)=0.$ From here and
(\ref{a53}) it follows that $G-c_{1}$log$k$ has a limit as $k\rightarrow 0$
and
\begin{equation}
\underset{k\rightarrow 0}{\lim }\left( G-c_{1}\text{log}k\right) =c_{1}\text{%
log}\left| \xi \right| +c_{2}+o(\left| \xi \right| ^{-1})\text{ \ \ as }%
\left| \xi \right| \rightarrow \infty .\text{ }  \label{a56}
\end{equation}
On the other hand, from (\ref{16}) it follows that the limit $%
\lim_{k\rightarrow 0}[G-u_{1}\left( 0,\xi \right) $log$k]$ exists and is
equal to $u_{2}\left( 0,\xi \right) .$ It can only happen if $u_{1}\left(
0,\xi \right) =c_{1}$ and $u_{2}\left( 0,\xi \right) $ is given by (\ref{a56}%
).

The proof is completed.

\textbf{III. Analytic properties of the resolvent of the Schr\"{o}dinger
operator.} Let $R_{\lambda }=\left( H-\lambda \right) ^{-1}$ be the
resolvent of the perturbed operator $H=\left( -\Delta +q\right) $, where the
support of $q$ belongs to a square
\begin{equation}
S=\{\xi \in Z^{2}:|\xi _{i}|\leq m\}.  \label{s}
\end{equation}
We are going to study the kernel $R_{k^{2}}(\xi ,\eta )$ of the operator $%
R_{k^{2}}$ (this kernel is Green's function for $H-k^{2}$) when both $\xi
,\eta \in S,$ and therefore we shall consider the truncated resolvent $%
\widehat{R}_{k^{2}}=\chi R_{k^{2}}\chi ,$ where $\chi $ is the
characteristic function of $S.$

Since the operator $H$ with a potential $q,$ whose support is bounded, may
have at most a finite number of eigenvalues, the resolvent $R_{\lambda }$
(and therefore, the operator $\widehat{R}_{\lambda })$ is meromorphic in $%
\lambda $ in the complex $\lambda $-plane with the cut along the segment $%
[0,8],$ and it has at most a finite number of poles which are eigenvalues of
$H$. The following Proposition is an immediate consequence of the statement
above.

\begin{proposition}
\label{p1}. The operators $R_{k^{2}}$ and $\widehat{R}_{k^{2}}$ are
meromorphic in $k\in \mathbf{C\backslash }\left[ -\sqrt{8},\sqrt{8}\right] .$
Outside of the segment $\left[ -\sqrt{8},\sqrt{8}\right] $ they have poles
of the first order at the points $k=\pm i\sigma _{j}$ for which $\lambda
_{j}=-\sigma _{j}^{2}$ are negative eigenvalues of $H$ and at points $k=\pm
\rho _{j}$ for which $\lambda _{j}=\rho _{j}^{2}>8$ are positive eigenvalues
of $H.$
\end{proposition}

The following statement was proved in \cite{s}:

\begin{proposition}
\label{p2}. The operator $H$ does not have eigenvalues on the set $(0,4)\cup
(4,8).$
\end{proposition}

It means that $H$ may have at most three eigenvalues imbedded into the
continuous spectrum, and those possible eigenvalues are $\lambda =0,4,8.$

We consider the upper half plane $\mathbf{C}_{+}=\{k:$Im$k>0\}$ as the image
of the set $\{\lambda :0<$arg$\lambda <2\pi \}$ under the map $k=\sqrt{%
\lambda }$ . Thus, the resolvent set of the operator $H$ in the $k$-plane is
$\overline{\mathbf{C}}_{+}\backslash \left[ -\sqrt{8},\sqrt{8}\right] $,
where the points $k=\pm \sigma $ with the same $\sigma >\sqrt{8}$ correspond
to the same value of $\lambda $. The next statement concerns an analytic
extension of the operator $\widehat{R}_{k^{2}}$ from the upper half plane $%
\mathbf{C}_{+}$ through the segment $\left[ -\sqrt{8},\sqrt{8}\right] .$ We
denote by $P_{j},$ $j=0,1,2,$ the orthogonal projections (in $L^{2}(Z^{2})$)
onto the eigenspaces of the operator $H$ with the eigenvalues $\lambda
=0,4,8,$ respectively. If some of those values of $\lambda $ are not
eigenvalues then the corresponding operators set to be zero. Let $\widehat{P}%
_{j}=\chi P_{j}\chi .$ Let $\widehat{R}_{k^{2}}($mod$P)$ be the operator
function $\widehat{R}_{k^{2}}$ by modulus of polynomial in $k$ operator
functions.

\begin{theorem}
\label{t21}1). The operator $\widehat{R}_{k^{2}},$ defined for $k\in
\overline{\mathbf{C}}_{+}\backslash \left[ -\sqrt{8},\sqrt{8}\right] ,$
admits an analytic extension on $\left[ -\sqrt{8},\sqrt{8}\right] \backslash
\cup k_{s}$, where $k_{0}=0,$ $k_{\pm 1}=\pm 2$, and $k_{\pm 2}=\pm \sqrt{8}%
. $

2). It has branch points of the logarithmic type at $k=k_{s},$ and there are
integers $\alpha _{s}$ and $\beta _{s}$ and operators $A_{s}\neq 0$ and $%
B_{s}$ ( in the finite dimensional space $L$ of functions supported on $S$)
such that the truncated resolvent $\widehat{R}_{k^{2}}$ has the following
behavior when $k\rightarrow k_{s}$, $k\neq k_{s}-i\sigma ,$ $\sigma \geq 0:$%
\begin{eqnarray}
\widehat{R}_{k^{2}}(\text{mod}P) &=&\widehat{P}_{|s|}\frac{1}{k_{s}^{2}-k^{2}%
}+A_{s}\left( k-k_{s}\right) ^{\alpha _{s}}\text{log}^{\beta _{s}}\left(
k-k_{s}\right) +B_{s}\left( k-k_{s}\right) ^{\alpha _{s}}\text{log}^{\beta
_{s}-1}\left( k-k_{s}\right)  \notag \\
&&+O\left( \left( k-k_{s}\right) ^{\alpha _{s}}\text{log}^{\beta
_{s}-2}\left( k-k_{s}\right) \right) \text{,}  \label{32}
\end{eqnarray}

3). The kernels $\widehat{R}_{k^{2}}(\xi ,\eta )$ of the operators $\widehat{%
R}_{k^{2}}$ with $k=\pm \sigma +i\varsigma ,$ where $\sigma $ is real, $%
\varsigma \geq 0,$ are complex adjoint (and therefore, $\alpha _{s}=\alpha
_{|s|},$ $\beta _{s}=\beta _{|s|})$.

The following inequalities hold: \ a) $\alpha _{0}\geq -2,$ and $\beta
_{0}\leq -1$ when $\alpha _{0}=-2;$ b) \ Let $s=\pm 1,\pm 2.$ Then $\alpha
_{s}\geq -1,$ and $\beta _{s}\leq -1$ when $\alpha _{s}=-1.$

4). If $q(\xi )\geq 0$ and $q$ is not equal to zero identically, then $%
\widehat{R}_{k^{2}}$ is bounded in a neighborhood of $k=0.$ In particular,
in this case $\lambda =0$ is not an eigenvalue of $H$ and $\alpha _{0}\geq
0. $
\end{theorem}

\textbf{Remark.} The reason to consider $\widehat{R}_{k^{2}}($mod$P)$ in (%
\ref{32}), but not $\widehat{R}_{k^{2}}$, is the following: the large time
asymptotic behavior of the solutions $v=v(t,x)$ to the Cauchy problem (\ref
{1}) depends on the first singular in $k-k_{s}$ term in the asymptotic
expansion for $\widehat{R}_{k^{2}}-\widehat{P}_{|s|}\frac{1}{k_{s}^{2}-k^{2}}
$ as $k\rightarrow k_{s}$. In fact, (\ref{32}) will be proven for $\widehat{R%
}_{k^{2}},$ but it will be clear from the proof that one can omit integer
nonnegative powers of $k-k_{s}$ and get (\ref{32}) for $\widehat{R}_{k^{2}}($%
mod$P).$

\textbf{Proof.} Let
\begin{equation}
\left( -\Delta -k^{2}+q\right) u=f,\text{ \ \ \ Im}k>0,  \label{35}
\end{equation}
where supports of $q$ and $f$ belong to $S$. The formula
\begin{equation}
u=R_{k^{2}}^{0}h,\text{ \ \ \ Im}k>0,  \label{for}
\end{equation}
establishes a one-to-one correspondence between solutions $u=R_{k^{2}}f\in
L^{2}(Z^{2})$ of (\ref{35}) and functions $h$ which are supported on $S$ and
satisfy the equation
\begin{equation}
h+qR_{k^{2}}^{0}h=f\ ,\text{ \ \ \ Im}k>0.\   \label{36}
\end{equation}
In fact, if $u\in L^{2}(Z^{2})$ is a solution of (\ref{35}), then $\left(
-\Delta -k^{2}\right) u=h$ with $h=f-qu.$ Hence, $u=R_{k^{2}}^{0}h$ and the
support of $h$ belongs to $S.$ The substitution of (\ref{for}) into (\ref{35}%
) leads to the equation (\ref{36}) for $h.$ On the other hand, if $h$ is
supported on $S$ and (\ref{36}) holds, then $\ u,$ given by (\ref{for}),
belongs to $L^{2}(Z^{2})$ and satisfies (\ref{35}). The latter can be
checked by substituting (\ref{for}) into (\ref{35}).

The equation (\ref{36}) has the form $T_{k}h=f$, where the operator $T_{k}$
acts in the finite dimensional space $L$ of functions supported on $S$. It
also can be written in the form of a linear system:
\begin{equation}
h(\xi )+\underset{\eta \in S}{\sum }q(\xi )G(k,\xi -\eta )h(\eta )=f(\xi ),%
\text{ \ \ \ }\xi \in S,\text{ \ \ \ Im}k>0,  \label{37}
\end{equation}
where $G$ is Green's function of the operator $-\Delta -k^{2}$ (the kernel
of $R_{k^{2}}^{0}).$ If one chooses the basis in $L,$ which consists of
functions equal to one at one point of $S$ and equal to zero everywhere
else, then the matrix $[T_{k}]$ of the operator $T_{k}$ in this basis is:
\begin{equation*}
\lbrack T_{k}]=[T_{k}(\xi ,\eta )]_{\xi ,\eta \in S}=E+[q(\xi )G(k,\xi -\eta
)],
\end{equation*}
where $E$ is the identity matrix. Note also that (\ref{for}) implies that
\begin{equation}
\widehat{R}_{k^{2}}=\chi R_{k^{2}}^{0}T_{k}^{-1},\text{ \ \ \ \ Im}k>0.
\label{121}
\end{equation}

From the first statement of Theorem \ref{t1} it follows that the matrix $%
[T_{k}]$ can be analytically extended on $\left[ -\sqrt{8},\sqrt{8}\right]
\backslash \cup k_{s}.$ The equation (\ref{35}) is uniquely solvable if $%
k\neq i\sigma _{j},$ and therefore the same is true for the equation (\ref
{36}). Thus, $\det [T_{k}]\neq 0$ when Im$k>0,$ $k\neq i\sigma _{j}.$ Hence,
the operator $T_{k}^{-1}$ is meromorphic on $\left[ -\sqrt{8},\sqrt{8}\right]
\backslash \cup k_{s}.$ From here, (\ref{121}), and the first statement of
Theorem \ref{t1} it follows that $\widehat{R}_{k^{2}}$ has a meromorphic
extension on $\left[ -\sqrt{8},\sqrt{8}\right] \backslash \cup k_{s}.$ On
the other hand, it was proved in (\cite{s}) that $\lim_{k\rightarrow
k^{0}+i0}\widehat{R}_{k^{2}}$ exists for any $k^{0}\in \left[ -\sqrt{8},%
\sqrt{8}\right] \backslash \cup k_{s}.$ This justifies the first statement
of Theorem \ref{t21}.

The second statement of Theorem \ref{t21} can be proven absolutely
similarly. The only difference is that the analytic extension of $[T_{k}]$
has branch points of the logarithmic type at $k=k_{s},$ $s=0,\pm 1,\pm 2$.
To be more exact, the elements of $[T_{k}]$ in neighborhoods of the points $%
k=k_{s}$ are linear functions of log$(k-k_{s})$ with coefficients analytic
in $k$ (see (\ref{16})). Let $N$ be the number of points in $S$. By solving
the system $[T_{k}]h=f$ using the Kramer rule we get that $T_{k}^{-1}$ in a
neighborhood of the point $k=k_{s}$ has the form
\begin{equation}
T_{k}^{-1}=\frac{\widehat{A}_{s}(k,\text{log}(k-k_{s}))}{B_{s}(k,\text{log}%
(k-k_{s}))},\text{ \ }s=0,\pm 1,\pm 2,  \label{38}
\end{equation}
where the functions $B_{s}=\det [T_{k}]$ are polynomials of order $N$ with
respect to the second argument (i.e., log$(k-k_{s}))$ whose coefficients are
analytic in $k$, and $\widehat{A}_{s}$ are polynomials in log$(k-k_{s})$ of
order $N-1$, whose coefficients are linear operators on $L$ (matrices) which
depend analytically in $k$. From here, (\ref{121}), and Theorem \ref{t1} it
follows that the formula (\ref{38}) is valid for $\widehat{R}_{k^{2}}$ with
the only difference that the order of the polynomials in the numerator is $%
N: $%
\begin{equation}
\widehat{R}_{k^{2}}=\frac{\widehat{C}_{s}(k,\text{log}(k-k_{s}))}{D_{s}(k,%
\text{log}(k-k_{s}))},  \label{bb}
\end{equation}
where the functions $D_{s}$ are polynomials in log$(k-k_{s})$ with analytic
in $k$ coefficients, and $\widehat{C}_{s}$ are polynomials in log$(k-k_{s})$
whose coefficients are analytic in $k$ operators. Then a similar
representation is valid for the difference of $\widehat{R}_{k^{2}}$ and $%
\widehat{P}_{|s|}\frac{1}{k_{s}^{2}-k^{2}}:$%
\begin{equation}
\widehat{R}_{k^{2}}-\widehat{P}_{|s|}\frac{1}{k_{s}^{2}-k^{2}}=\frac{%
\widehat{C}_{s}}{D_{s}}-\frac{\widehat{P}_{|s|}}{k_{s}^{2}-k^{2}}=\frac{%
\widehat{F}_{s}(k,\text{log}(k-k_{s}))}{G_{s}(k,\text{log}(k-k_{s}))},
\label{pln}
\end{equation}
where $\widehat{F}_{s}=\left( k_{s}^{2}-k^{2}\right) \widehat{C}_{s}-D_{s}%
\widehat{P}_{|s|},$ $G_{s}=\left( k_{s}^{2}-k^{2}\right) D_{s}.$ The formula
(\ref{32}) immediately follows from here. The second statement of Theorem
\ref{t21} is proved.

The first part of the third statement of the Theorem is an obvious
consequence of the fact that
\begin{equation}
R_{\overline{\lambda }}=\overline{R_{\lambda }}\text{ \ \ if \ }\lambda
\notin \lbrack 0,\infty ).  \label{mmn}
\end{equation}
Now we are going to prove the second part. Let $dE_{\sigma },$ $\sigma $ is
real, be the operator valued spectral measure for the operator $H$. Let $\mu
_{ac}(d\sigma )$ be the absolutely continuous part of the measure $%
dE_{\sigma },$ and let the operator valued function $\nu =\nu (\sigma )$ be
the density of the measure $\mu _{ac}(d\sigma ).$ The matrix elements $\nu
_{\xi ,\eta },$ $\xi ,\eta \in Z^{2},$\ of the operator $\nu $ are the
densities of the absolutely continuous parts of the scalar measures $%
(dE_{\sigma }\delta _{\xi },\delta _{\eta }),$ where $\delta _{\xi }$ is the
function on $Z^{2}$ equal to one at a fixed point $\xi $ and equal zero
elsewhere. We also shall consider the restriction $\widehat{\nu }=\chi \nu
(\sigma )\chi $ of the operators $\nu =\nu (\sigma )$ onto the space $L$ of
functions $\psi $ supported on $S.$ From Stone's formula and the first part
of Theorem \ref{t21}\ it follows that the density $\widehat{\nu }(\sigma )$
in a neighborhood of $\sigma =\lambda _{s}$ can be obtained as
\begin{equation}
\widehat{\nu }(\sigma )=\frac{1}{\pi }\lim_{\lambda \rightarrow \sigma +i0}%
\text{Im}[\widehat{R}_{\lambda }-\widehat{P}_{|s|}\frac{1}{\lambda
_{s}-\lambda }].  \label{vu1}
\end{equation}
Obviously, the operator function $\nu $ is summable in $\sigma $:
\begin{equation}
\widehat{\nu }\in L^{1}.  \label{c2}
\end{equation}

The formula (\ref{pln}) implies a more exact asymptotic expansion of $%
\widehat{R}_{k^{2}\text{ }}$ than the expansion (\ref{32}). Namely,
\begin{equation}
\widehat{R}_{k^{2}}-\widehat{P}_{|s|}\frac{1}{k_{s}^{2}-k^{2}}\text{ }\sim
\text{ }\sum_{j\geq 0}(k-k_{s})^{\alpha _{s}+j}\frac{\widehat{T}_{j,s}(\text{%
log}(k-k_{s}))}{N_{j,s}(\text{log}(k-k_{s}))},\text{ \ \ \ }k\rightarrow
k_{s},  \label{kon}
\end{equation}
where $N_{j,s}$ are polynomials, $\widehat{T}_{j,s}$ are operator valued
polynomials, and $\widehat{T}_{0,s}$ is a non-zero operator function.
Consider first the case $s=1$ or $2$. Then (\ref{kon}) implies that the
following expansion is valid when $\lambda \ $is close to $\lambda _{s}$ and
Im$\lambda >0:$
\begin{equation}
\widehat{R}_{\lambda }-\widehat{P}_{|s|}\frac{1}{\lambda _{s}-\lambda }\text{
}\sim \text{ }\sum_{j\geq 0}(\lambda -\lambda _{s})^{\alpha _{s}+j}\frac{%
\widehat{T}_{j,s}(\text{log}(\lambda -\lambda _{s}))}{N_{j,s}(\text{log}%
(\lambda -\lambda _{s}))},  \label{3221}
\end{equation}
where the polynomials $\widehat{T}_{j,s},$ $N_{j,s}$ are different from
those in (\ref{kon}), and still $\widehat{T}_{0,s}\neq 0$. \ Together with (%
\ref{vu1}) this leads to the following representation
\begin{equation}
\widehat{\nu }(\sigma )=\lim_{\lambda \rightarrow \sigma +i0}\text{Im\{}%
(\lambda -\lambda _{s})^{\alpha _{s}}\text{log}^{\beta _{s}}(\lambda
-\lambda _{s})f(\text{log}^{-1}(\lambda -\lambda _{s}))+O((\lambda -\lambda
_{s})^{\alpha _{s}+1}\text{log}^{\gamma _{s}}(\lambda -\lambda _{s}))\text{%
\},}  \label{c1}
\end{equation}
where the operator function $f$ is analytic, $f(0)\neq 0$, and $\gamma _{s}$
is an integer.

Assume that $\alpha _{s}\leq -2,$ $\widetilde{\beta }_{s}\neq 0.$ From (\ref
{c1}) with $\sigma >\lambda _{s}$ and (\ref{c2}) we obtain that the
operators $f(\sigma )$ with real $\sigma >\lambda _{s}$ are real valued
(i.e. they map any real valued function into a real valued function). In
particular, the operators $f^{(j)}(0)$ are real valued for all $j\geq 0$.
Then (\ref{c1}) with $\sigma <\lambda _{s}$ implies that
\begin{eqnarray}
\widehat{\nu }(\sigma ) &=&(\sigma -\lambda _{s})^{\alpha _{s}}\text{Im}[%
\text{log}^{\beta _{s}}(|\sigma -\lambda _{s}|+i\pi )f(0)+\text{log}^{\beta
_{s}-1}(|\sigma -\lambda _{s}|+i\pi )f^{\prime }(0)  \notag \\
&&+O(\text{log}^{\beta _{s}-2}(|\sigma -\lambda _{s}|)]  \notag \\
&=&\pi (\sigma -\lambda _{s})^{\alpha _{s}}[\beta _{s}f(0)\text{log}^{\beta
_{s}-1}|\sigma -\lambda _{s}|+O(\text{log}^{\beta _{s}-2}(|\sigma -\lambda
_{s}|)],\text{ \ \ }\sigma <\lambda _{s}.  \label{c3}
\end{eqnarray}
This contradicts (\ref{c2}), and therefore the assumption is wrong. Let us
assume that $\alpha _{s}\leq -2,$ $\beta _{s}=0.$ Let $n$ be the smallest $%
j>0$ for which $f^{(j)}(0)\neq 0.$ The formula (\ref{c1}) leads to the
following analog of (\ref{c3}):
\begin{equation*}
\widehat{\nu }(\sigma )=\pi (\sigma -\lambda _{s})^{\alpha _{s}}[\frac{-n}{n!%
}f^{(n)}(0)\text{log}^{-n-1}|\sigma -\lambda _{s}|+O(\text{log}%
^{-n-2}(|\sigma -\lambda _{s}|)],\text{ \ \ }\sigma <\lambda _{s},
\end{equation*}
which also contradicts (\ref{c2}). Thus $\alpha _{s}\leq -2$ requires $\beta
_{s}=0,$ $f$ is a real valued operator which does not depend on $\lambda $.
Now the same arguments can be applied successively to the next terms in (\ref
{3221}) for which $\alpha _{s}+j<-1$. When $\alpha _{s}+j=-1,$ we can only
get that
\begin{equation}
\frac{\widehat{T}_{j,s}}{N_{j,s}}\text{ }=\text{ }M+O(\text{log}%
^{-1}(\lambda -\lambda _{s})),\text{ \ \ }\alpha _{s}+j=-1,  \label{dd}
\end{equation}
where $M$ is a real valued operator. So, we obtain that (\ref{3221}) has the
form
\begin{equation}
\widehat{R}_{\lambda }-\widehat{P}_{|s|}\frac{1}{\lambda _{s}-\lambda }\text{
}=\text{ }\sum_{-\alpha _{s}>j\geq 0}(\lambda -\lambda _{s})^{\alpha
_{s}+j}M_{j,s}+O((\lambda -\lambda _{s})^{-1}\text{log}^{-1}(\lambda
-\lambda _{s}))  \label{dd1}
\end{equation}
where $M_{j,s}$ are real valued and independent of $\lambda $ operators on $%
L $. If $\alpha _{s}=-1,$ then the arguments above still lead to (\ref{dd})
when $j=0.$ Thus (\ref{dd1}) is valid when $\alpha _{s}\leq -1.$

Our next goal is to show that (\ref{dd1}) with $\alpha _{s}<-1$ implies that
$M_{0,s}$ are zero operators. Let $\delta >0$ be so small that the interval $%
\omega =[\lambda _{s},\lambda _{s}+\delta ]$ does not contain the
eigenvalues of $H$ different from $\lambda _{s}.$ Let $\omega _{\varepsilon
}=$ $[\lambda _{s}+i\varepsilon ,\lambda _{s}+\delta +i\varepsilon ]$ be the
shift of $\omega .$ Stone's formula implies that
\begin{equation}
\int_{\omega }\widehat{\nu }(\sigma )d\sigma =\frac{1}{\pi }%
\lim_{\varepsilon \rightarrow +0}\text{Im}\int_{\omega _{\varepsilon }}[%
\widehat{R}_{\lambda }-\widehat{P}_{|s|}\frac{1}{\lambda _{s}-\lambda }%
]d\lambda .  \label{as}
\end{equation}
Assume that $M_{0,s}\neq 0.$ Then the existence of the limit in the right
hand side above contradicts (\ref{dd1}) if $\alpha _{s}<-1$ is odd. It also
contradicts (\ref{dd1}) if $\alpha _{s}<-1$ is even and $M_{1,s}\neq 0.$ One
can easily show that $\alpha _{s}$ does not depend on the choice of the
square $S$ if $S$ contains the support of $q$ (this follows from (\ref{121})
and (\ref{37})). Hence, one can take a bigger $S$, so that there is a
function $\psi $ for which $M_{0,s}\psi $ is not identically equal to zero
on a smaller cube $S^{\prime }=\{\xi :|\xi _{i}|\leq m-1\}\subset S$ (see (%
\ref{s})). We apply the operator
\begin{equation*}
H-\lambda =-\Delta +q(\xi )-\lambda _{s}-(\lambda -\lambda _{s})
\end{equation*}
to both sides of (\ref{dd1}). The following relations hold in $S^{\prime }$
for any function $\psi $ defined on $Z^{2}:$
\begin{equation*}
(H-\lambda )\widehat{R}_{\lambda }\psi =\psi ,\text{ \ }(H-\lambda )\widehat{%
P}_{|s|}\frac{\psi }{\lambda _{s}-\lambda }=-\widehat{P}_{|s|}\psi
,~~\rm{ Im}\lambda >0.
\end{equation*}
The right hand sides here are independent of $\lambda .$ Thus, the negative
powers of $\lambda -\lambda _{s}$ will be cancelled when $H-\lambda $ is
applied to the right hand side of (\ref{dd1}). \ This leads to the following
relation in $S^{\prime }$:
\begin{equation*}
(H-\lambda _{s})M_{1,s}\psi =M_{0,s}\psi ,\text{ \ if }\alpha _{s}\leq -2.
\end{equation*}
Hence, the assumption $M_{0,s}\neq 0$ implies that $M_{1,s}\neq 0,$\ and
this leads to the contradiction discussed above. This proves that $\alpha
_{s}\geq -1.$

If $\alpha _{s}=-1,$ then (\ref{dd1}) implies that the limit in the right
hand side of (\ref{as}) exists and converges to $\frac{1}{2}M_{0,s}$ when $%
\delta =|\omega |\rightarrow 0.$ The left hand side in (\ref{as}) converges
to zero as $\delta \rightarrow 0.$ Thus, $M_{0,s}=0$ and (\ref{dd1}) takes
the form
\begin{equation*}
\widehat{R}_{\lambda }-\widehat{P}_{|s|}\frac{1}{\lambda _{s}-\lambda }\text{
}=\text{ }O((\lambda -\lambda _{s})^{-1}\text{log}^{-1}(\lambda -\lambda
_{s})),\text{ \ }s=1,2.
\end{equation*}
This proves the second part of statement 3 of Theorem \ref{t21} in the case $%
s=1,2.$ If $s=-1$ or $-2$, the result follows from (\ref{mmn}). The
statement for $s=0$ can be justified similarly to the case $s=1,2.$ We leave
it for the reader. This completes the proof of statement 3.

Next we shall show that $\widehat{R}_{k^{2}}$ is bounded in a neighborhood
of $k=0$ when $q\geq 0$. From (\ref{bb}) it follows that
\begin{equation}
\widehat{R}_{k^{2}}=Ak^{\alpha }\text{log}^{\beta }k+Bk^{\alpha }\text{log}%
^{\beta -1}k+O(k^{\alpha }\text{log}^{\beta -2}k)\text{ \ as }k\rightarrow 0.
\label{a58}
\end{equation}
where $A$ is a non zero operator defined on the space $L$\ of functions with
supports in $S$. One could get (\ref{a58}) from (\ref{32}) with $\alpha
=\alpha _{0},$ $\beta =\beta _{0}$ (when $P_{0}=0$)$,$ or $\alpha =-2,$ $%
\beta =0$ (if $P_{0}\neq 0$)$.$ It is also possible to get (\ref{a58}) from (%
\ref{bb}). The last statement of the Theorem will be proved if we show that
the first term in the right hand side of (\ref{a58}) is bounded as $%
k\rightarrow 0.$

Suppose, to the contrary, that $\alpha <0$ or $\alpha =0,$ $\beta >0$. Let $%
f $ be a function supported on $S$, such $Af$ is a non-zero function. We
have
\begin{equation}
\left( -\Delta -k^{2}\right) R_{k^{2}}f=f-qR_{k^{2}}f,\text{ \ \ Im}k>0.
\label{equ}
\end{equation}
Since $q$ and $f$ are supported on $S$, one can replace $R_{k^{2}}$ in the
right hand side above by $\widehat{R}_{k^{2}}.$ Thus,
\begin{equation*}
R_{k^{2}}f=R_{k^{2}}^{0}(f-q\widehat{R}_{k^{2}}f),\text{ \ \ Im}k>0.
\end{equation*}
From here, (\ref{a58}), and Lemma \ref{l00} it follows that, for any $\xi
\in Z^{2},$ Im$k>0$ and $k\rightarrow 0,$ the function $(R_{k^{2}}f)(\xi )$
has the form
\begin{equation}
R_{k^{2}}f=[c_{1}\text{log}k+u_{2}(0,\xi )+O(k\text{log}k)]\ast
\{f-q[Afk^{\alpha }\text{log}^{\beta }k+Bfk^{\alpha }\text{log}^{\beta
-1}k+O\left( k^{\alpha }log^{\beta -2}k\right) ]\},  \label{tr}
\end{equation}
where the convolution of two functions $g,$ $h$ on $Z^{2}$ is defined as
\begin{equation*}
(g\ast h)(\xi )=\sum_{\eta }g(\xi -\eta )h(\eta ).
\end{equation*}
Hence, the main term in the asymptotic expansion of $R_{k^{2}}f$ as $%
k\rightarrow 0$ is
\begin{equation}
-c_{1}k^{\alpha }\text{log}^{\beta +1}k\sum_{\eta }(qAf)(\eta ).  \label{qw}
\end{equation}
On the other hand, $(R_{k^{2}}f)(\xi )=(\widehat{R}_{k^{2}}f)(\xi )$ when $%
\xi \in S,$ and this value has a smaller order due to (\ref{a58}). Thus (\ref
{qw}) is zero. Since $c_{1}\neq 0,$ we arrive at
\begin{equation}
\sum_{\eta }(qAf)(\eta )=0.  \label{in}
\end{equation}

Now from (\ref{tr}) we get
\begin{equation}
R_{k^{2}}f=v(\xi )k^{\alpha }\text{log}^{\beta }k+O(k^{\alpha }\text{log}%
^{\beta -1}k),\text{ \ \ Im}k>0,\text{ \ }k\rightarrow 0,  \label{rrr}
\end{equation}
where
\begin{equation*}
v(\xi )=-u_{2}(0,\xi )\ast qAf-c_{1}\ast qBf
\end{equation*}
Together with Lemma \ref{l00} this leads to the following asymptotics for $%
v(\xi ):$
\begin{equation}
v(\xi )=-c_{1}\ln |\xi |\ast qAf+c_{2}+o(|\xi |^{-1}),\text{ \ }|\xi
|\rightarrow \infty .  \label{vv}
\end{equation}

We will need more specific behavior of $v$ at infinity. First of all note
that (\ref{in}) implies that
\begin{equation*}
\ln |\xi |\ast qAf=\sum_{\eta }(\ln |\xi -\eta |-\ln |\xi |)qAf(\eta
)=O(|\xi |^{-1}),\text{ \ }|\xi |\rightarrow \infty ,
\end{equation*}
and therefore,
\begin{equation}
|v(\xi )|<C<\infty .  \label{v}
\end{equation}
Let us denote by $\frac{\partial }{de},$ $e\in Z^{2},$ the difference
derivative in the direction of $e:$%
\begin{equation}
\frac{\partial h(\xi )}{de}=h(\xi +e)-h(\xi ).  \label{der}
\end{equation}
Let us show that
\begin{equation}
\frac{\partial v(\xi )}{de}=o(|\xi |^{-1}),\text{ \ }|\xi |\rightarrow
\infty .  \label{dv}
\end{equation}
Formula (\ref{in}) implies that
\begin{eqnarray*}
\frac{\partial }{de}\ln |\xi |\ast qAf &=&\sum_{\eta }(\ln |\xi +e-\eta
|-\ln |\xi -\eta |)qAf(\eta ) \\
&=&\sum_{\eta }(\ln |\xi +e-\eta |-\ln |\xi -\eta |-\ln |\xi +e|+\ln |\xi
|)qAf(\eta )=O(|\xi |^{-2}).
\end{eqnarray*}
The last equality can be easily obtained by evaluating the asymptotics for
the combination of the logarithms above. The formula above and (\ref{vv})
justify (\ref{dv}).

Now we multiply both sides of (\ref{equ}) by $k^{-\alpha }log^{-\beta }k$
and pass to the limit as $k\rightarrow 0.$ This leads to
\begin{equation}
\left( -\Delta +q\right) v=0\,.  \label{a61}
\end{equation}

We apply the Green's Theorem for the difference Schr\"{o}dinger operator to $%
v$ over the square $-N\leq \xi _{1},\xi _{2}\leq N$. Estimates (\ref{v}), (%
\ref{vv}) allow us to pass to the limit as $N\rightarrow \infty $. This
leads to the following result:\
\begin{equation}
0=\sum_{\xi \in Z^{2}}\left( -\Delta +q\right) v\cdot \overline{v}=\sum_{\xi
\in Z^{2}}(|\nabla _{\xi }v|^{2}+q|v|^{2}),  \label{a62}
\end{equation}
where $\nabla _{\xi }=(\frac{\partial }{\partial e_{1}},\frac{\partial }{%
\partial e_{2}}),$ $e_{1}=(1,0),$ $e_{2}=(0,1)$[(+,)] and $\frac{\partial }{%
\partial e}$ is defined in (\ref{der}). Since $q\geq 0$, (\ref{a62}) implies
that $v$ is a constant. Now from (\ref{a61}) (or (\ref{a62})) it follows
that $v=0,$ since $q$ is a non-zero function. This conclusion contradicts (%
\ref{rrr}) and (\ref{a58}). In fact, the later two relations imply that $%
v(\xi )=Af$ when $\xi \in S.$ Yet, $f$ was chosen so that $Af$ \ is a
non-zero function on $S$. The contradiction proves that our assumption of
unboundedness of the first term in the right hand side of (\ref{a58}) is
wrong. This completes the proof of the last statement of Theorem \ref{t21}.

\textbf{IV. The large time behavior of solutions to the Cauchy problem. }The
solution $v=v\left( t,\xi \right) $ to the Cauchy problem (\ref{1}) is given
by (\ref{9}). We assume that $f$ and $q$ are real valued and have bounded
supports, and we will study the asymptotic behavior of $v$ when $%
t\rightarrow \infty $ and $\xi \in S$, where $S$ is a square in $Z^{2}$
which contains the supports of $q$ and $f.$ In order to keep in mind this
restriction on $\xi $ we multiply both sides of \ (\ref{9}) by the
characteristic function $\chi $ of $S$. Since the support of $f$ belongs to $%
S$, we can replace $f$ in (\ref{9}) by $\chi f$, and we arrive to
\begin{equation}
\chi v\left( t,\xi \right) =\frac{1}{2\pi }\int_{B-i\infty }^{B+i\infty }%
\widehat{R}_{k^{2}}fe^{-ikt}dk,\text{ \ \ }B>A.  \label{q1}
\end{equation}

Let us denote by $P_{(j)},$ $1\leq j\leq n_{1},$ the projection operators in
$L^{2}(Z^{2})$ onto the eigenspaces of the operator $H$ with the negative
eigenvalues $\lambda =-\sigma _{j}^{2}.$ Let $P_{j},$ $1\leq j\leq n_{2},$
be the projection operators onto the eigenspaces of the operator $H$ with
the positive eigenvalues $\lambda =\rho _{j}^{2}.$ We consider here all
positive eigenvalues, including those which belong to the continuous
spectrum (see the Proposition \ref{p2}$)$. Let $P_{0}$ be the projection
operator onto the eigenspace with the eigenvalue $\lambda =0$, and let $%
P_{0}=0$ if $\lambda =0$ is not an eigenvalue. Let $\widehat{P}_{(j)}=\chi
P_{(j)}\chi ,$ $\widehat{P}_{j}=\chi P_{j}\chi ,$ $\widehat{P}_{0}=\chi
P_{0}\chi .$

\begin{theorem}
\label{tlast}. The solution $v$ of the problem (\ref{1}) has the following
asymptotic behavior as $t\rightarrow \infty :$%
\begin{eqnarray*}
\chi v\left( t,\xi \right) &=&\sum_{j}\frac{e^{\sigma _{j}t}}{2\sigma _{j}}%
\widehat{P}_{(j)}f+\sum_{j}\frac{\sin \rho _{j}t}{\rho _{j}}\widehat{P}%
_{j}f+t\widehat{P}_{0}f \\
&&+\sum_{0\leq s\leq 2}a_{s}t^{-\alpha _{s}-1}\left( \text{log}t\right)
^{\beta _{s}-\gamma _{s}}\sin (k_{s}t+\omega _{s})A^{s}f+w\left( t,\xi
\right) ,
\end{eqnarray*}
where $a_{s}\neq 0,$ $\omega _{s}$ are some constants, $\alpha _{s},$ $\beta
_{s},$ $A^{s}$ are defined in (\ref{32}), $\gamma _{s}=0$ if $\alpha
_{s}=-1, $ $\gamma _{s}=1$ if $\alpha _{s}\geq 0,$ $k_{0}=0,$ $k_{1}=2,$ $%
k_{2}=\sqrt{8},$ and $\frac{d^{n}}{dt^{n}}w(t,\xi )$ decay as $t\rightarrow
\infty $ by the factor log$^{-1}t$ faster than the slowest term in the last
sum in the right hand side above.
\end{theorem}

\textbf{Proof. }We plan to move the contour of integration in (\ref{q1})
down and replace it by the following contour $\Gamma $. Let $c_{s}$ be the
following rays in the complex $k$-plane $\mathbf{C:}$%
\begin{equation*}
c_{s}=\{k=k_{s}-i\sigma ,\text{ \ }\sigma \geq 0\},\text{ \ \ }s=0,\pm 1,\pm
2,
\end{equation*}
and let $\mathbf{C}^{\prime }=$ $\mathbf{C\backslash \cup }c_{s}$ be the
complex $k$-plane $\mathbf{C}$ with the cuts along $c_{s}.$ Let $\Gamma _{0}$
be the line Im$k=-\delta ,$ $\delta >0,$ without points $k_{s}-i\delta $,
where the line intersects $c_{s}.$ Let $\gamma _{s}$ be smooth loops joining
the points $k_{s}-i\delta -0$ (on the left side of $c_{s}$) and $%
k_{s}-i\delta +0$ (on the right side of $c_{s}$) going round the cuts $c_{s}$
and lying in the $\delta $-neighborhood of the points $k_{s}.$ Then $\Gamma
=\Gamma _{0}\cup \gamma _{s}\ $is the union of $\Gamma _{0}\mathbf{\ }$and
all $\gamma _{s},$ and we need only to chose $\delta $ in order to specify $%
\Gamma .$ $\ $

From Proposition \ref{p1} it follows that there is a $\delta _{1}>0$ such
that the operator $\widehat{R}_{k^{2}}$ defined on $\mathbf{C}\backslash
\lbrack \sqrt{8},\sqrt{8}]$ does not have poles in the strip $-\delta
_{1}\leq $Im$k<0.$ Below we shall always consider $\widehat{R}_{k^{2}}$ in $%
\mathbf{C}^{\prime },$ i.e. we start with $\widehat{R}_{k^{2}}$ in the upper
half plane $\mathbf{C}_{+}$ and then extend it analytically down. Due to
Proposition \ref{p1}, this extension exists in the region $|$Re$k|>\sqrt{8}$
(which does not contain cuts $c_{s}$). According to Theorem \ref{t21}, there
exists $\delta _{2}>0$ such that $\widehat{R}_{k^{2}}$ has an analytic
(without poles) extension onto the $\delta _{2}$-neighborhood in $\mathbf{C}%
^{\prime }$ of the segment $[\sqrt{8},\sqrt{8}]$ with branch points at $%
k=k_{s}$. In fact, one can show that $\widehat{R}_{k^{2}}$ can be
meromorphically extended onto the whole $\mathbf{C}^{\prime },$ but we don't
need this statement. We choose an arbitrary $\delta <\min (\delta
_{1},\delta _{2}).$ Then $\widehat{R}_{k^{2}}$ is meromorphic strictly
inside the region $\Omega \subset \mathbf{C}^{\prime }$ between the contour
of integration in (\ref{q1}) and $\Gamma ,$and $\widehat{R}_{k^{2}}$ has
there poles only at points $k=i\sigma _{j},$ $\sigma _{j}>0$, and $k=\pm
\rho _{j}$ for which $\lambda =k^{2}$ are eigenvalues of $H$, which do not
belong to the absolutely continuous spectrum of $H$.

Note that
\begin{equation*}
||\widehat{R}_{k^{2}}||\leq ||R_{k^{2}}||\leq 1/d,\text{ \ \ }k\in \mathbf{C}%
\backslash \lbrack \sqrt{8},\sqrt{8}],
\end{equation*}
where $d$ is the distance in the complex $\lambda $-plane between the point $%
\lambda =k^{2}$ and the spectrum of the operator $H$. \ Thus,
\begin{equation}
||R_{k^{2}}||\leq C/|k|^{2},\text{ \ \ }k\in \Omega ,\text{ \ }%
|k|\rightarrow \infty .  \label{bes}
\end{equation}
Since $||\widehat{R}_{k^{2}}||\leq ||R_{k^{2}}||,$ the estimate (\ref{bes})
is valid also for $||\widehat{R}_{k^{2}}||.$ This allows us to use the
Cauchy theorem and reduce the integral (\ref{q1}) to the following form:
\begin{equation}
\chi v\left( t,\xi \right) =-i\sum_{j}\underset{k=i\sigma _{j}}{\text{Res}}%
\text{(}\widehat{R}_{k^{2}}fe^{-ikt})-i\widetilde{\sum_{j}}\underset{k=\pm
\rho _{j}}{\text{Res}}\text{(}\widehat{R}_{k^{2}}fe^{-ikt})+\frac{1}{2\pi }%
\int_{\Gamma }\widehat{R}_{k^{2}}fe^{-ikt}dk,  \label{res}
\end{equation}
where the sum does not include points $\pm \rho _{j}\in \lbrack -\sqrt{8},%
\sqrt{8}].$

Let us estimate the last term in (\ref{res}). We split $\Gamma _{0}=\Gamma
_{1}\cup \Gamma _{2},$ where $\Gamma _{1}$ is the bounded part of $\Gamma
_{0}$ for which $|$Re$k|<\sqrt{8}+1$ and $\Gamma _{2}$ is the part of $%
\Gamma _{0}$ where $|$Re$k|>\sqrt{8}+1.$ We denote
\begin{equation}
u_{s}\left( t,\xi \right) =\frac{1}{2\pi }\int_{\Gamma s}\widehat{R}%
_{k^{2}}fe^{-ikt}dk,\text{ \ \ }s=1,2,  \label{w2}
\end{equation}
Since $||\widehat{R}_{k^{2}}||$ is bounded and $|e^{-ikt}|=e^{-\delta t}$ on
$\Gamma _{1},$
\begin{equation*}
|\frac{\partial ^{j}}{\partial ^{j}t}u_{1}\left( t,\xi \right) |\leq \frac{1%
}{2\pi }\int_{\Gamma s}|k^{j}\widehat{R}_{k^{2}}fe^{-ikt}|dk\leq
C_{j}e^{-\delta t},\text{ \ \ \ }j=0,1,...\text{ .}
\end{equation*}
A similar estimate for $u_{2}$ with $j=0$ immediately follows from (\ref{bes}%
). To get the estimate for the derivatives we differentiate the equation
\begin{equation*}
\left( -\Delta +q+k^{2}\right) u=f,\text{ \ \ \ \ }k\in \Gamma _{2},\text{ \
\ }u=R_{k^{2}}f,\text{\ }
\end{equation*}
with respect to $k$:
\begin{equation*}
\left( -\Delta +q+k^{2}\right) u_{k}=-2ku.
\end{equation*}
From here and (\ref{bes}) \ it follows that $||\frac{\partial }{\partial k}%
R_{k^{2}}||\leq C/|k|^{3},$ \ \ $k\in \Gamma _{2}.$ If we differentiate the
equation again we can estimate the second derivative of the resolvent, and
so on. Thus,
\begin{equation}
||\frac{\partial ^{j}}{\partial k^{j}}R_{k^{2}}||\leq C/|k|^{2+j},\ \ k\in
\Gamma _{2}.  \label{ee}
\end{equation}
We apply the integration by parts $j$ times to the integral (\ref{w2}) with $%
s=2:$%
\begin{equation*}
u_{2}\left( t,\xi \right) =\sum_{1\leq l\leq j}[\frac{\partial ^{l-1}}{%
\partial k^{l-1}}(\widehat{R}_{k^{2}}f)\frac{e^{-ikt}}{(it)^{l}}]\left| _{k=-%
\sqrt{8}-1-i\delta }^{k=\sqrt{8}+1-i\delta }\right. +\frac{1}{2\pi }%
\int_{\Gamma _{2}}\frac{\partial ^{j}}{\partial k^{j}}(\widehat{R}_{k^{2}}f)%
\frac{e^{-ikt}}{(it)^{j}}dk.
\end{equation*}
This formula and (\ref{ee}) allow one to estimate the derivatives of $u_{2}:$%
\begin{equation*}
|\frac{\partial ^{j}}{\partial ^{j}t}w_{2}\left( t,\xi \right) |\leq
C_{j}e^{-\delta t},\text{ \ \ \ }j=0,1,...\text{ .}
\end{equation*}

Now (\ref{res}) takes the form
\begin{equation}
\chi v\left( t,\xi \right) =-i\sum_{j}\underset{k=i\sigma _{j}}{\text{Res}}%
\text{(}\widehat{R}_{k^{2}}fe^{-ikt})-i\widetilde{\sum_{j}}\underset{k=\pm
\rho _{j}}{\text{Res}}\text{(}\widehat{R}_{k^{2}}fe^{-ikt})+\sum_{s}\frac{1}{%
2\pi }\int_{\gamma _{s}}\widehat{R}_{k^{2}}fe^{-ikt}dk+u,  \label{lit}
\end{equation}
where $u=u_{1}+u_{2},$ and therefore
\begin{equation*}
|\frac{\partial ^{j}}{\partial ^{j}t}u\left( t,\xi \right) |\leq
C_{j}e^{-\delta t},\text{ \ \ \ }j=0,1,...\text{ .}
\end{equation*}

Since the operator $R_{\lambda }$ in a neighborhood of the eigenvalue $%
\lambda =\rho _{j}^{2}>8$ has the form
\begin{equation*}
R_{\lambda }=P_{j}\frac{1}{\rho _{j}^{2}-\lambda }+Q_{\lambda },
\end{equation*}
where $P_{j}$ is the projection operator onto the corresponding eigenspace
and $Q_{\lambda }$ is analytic in a neighborhood of $\rho _{j}^{2},$ we
obtain
\begin{equation*}
\underset{k=\rho _{j}}{\text{Res}}\text{(}\widehat{R}_{k^{2}}fe^{-ikt})+%
\underset{k=-\rho _{j}}{\text{Res}}\text{(}\widehat{R}_{k^{2}}fe^{-ikt})=%
\frac{i\sin \rho _{j}t}{\rho _{j}}\widehat{P}_{j}f\text{, \ \ \ \ where }%
\widehat{P}_{j}=\chi P_{j}\chi .
\end{equation*}
Similarly,
\begin{equation*}
\underset{k=i\sigma _{j}}{\text{Res}}\text{(}\widehat{R}_{k^{2}}fe^{-ikt})=%
\frac{ie^{\sigma _{j}t}}{2\sigma _{j}}\widehat{P}_{(j)}f.
\end{equation*}
Since the integral of any analytic (in particular, polynomial) function over
$\gamma _{s}$ is zero, we can replace $\widehat{R}_{k^{2}}$ in the last term
in (\ref{lit}) by the right hand side in (\ref{32}). After that the
corresponding integrals involving the first term from the right hand side of
(\ref{32}) can be evaluated with the help of the residue theorem. Together
with the two formulas above this allows us to rewrite (\ref{lit}) in the
form
\begin{equation}
\chi v\left( t,\xi \right) =\sum_{j=1}^{n_{1}}\frac{e^{\sigma _{j}t}}{%
2\sigma _{j}}\widehat{P}_{(j)}f+\sum_{j=1}^{n_{2}}\frac{\sin \rho _{j}t}{%
\rho _{j}}\widehat{P}_{j}f+t\widehat{P}_{0}f+\sum_{s}\frac{1}{2\pi }%
\int_{\gamma _{s}}Q_{s}(k)fe^{-ikt}dk+u,  \label{lll}
\end{equation}
where
\begin{equation*}
Q_{s}(k)=\widehat{R}_{k^{2}}(\text{mod}P)-\widehat{P}_{|s|}\frac{1}{%
k_{s}^{2}-k^{2}}
\end{equation*}
is defined in (\ref{32}), and the first three terms in the right hand side
of (\ref{lll}) contain the projections on all eigenspaces of $H$, including
those for which $\lambda _{j}\in \lbrack 0,8]$. In order to complete the
proof of the theorem, it remains only to evaluate the last term in (\ref{lll}%
) using Theorem \ref{t21} and the following simple statement which can be
found in \cite{v}, Ch X, Lemmas 7, 8:

\begin{lemma}
Let
\begin{equation*}
v(t)=\int_{l}\omega (k)e^{-ikt}dk,
\end{equation*}
where the function $\omega $ is analytic when $|k|<2\delta $ with a branch
point of the logarithmic type at $k=0,$ and $l$ is a smooth contour in the $%
\delta -$neighborhood of $k=0$ starting at $k=-i\delta -0,$ ending at $%
k=-i\delta +0$ and not having common points with the negative imaginary
semiaxis except the end points $k=-i\delta \pm 0.$ Let $\omega (k)$ have the
following asymptotic behavior as $k\rightarrow 0,$ $\frac{3\pi }{2}<$arg$k<-%
\frac{\pi }{2}:$%
\begin{equation*}
\omega (k)=k^{p}\text{log}^{q}k+ck^{p}\text{log}^{q-1}k+O(k^{p}log^{q-2}k)
\end{equation*}
where $p$ and $q$ are integers and $q\neq 0$ if $p\geq 0.$

Then
\begin{equation*}
v\left( t\right) =at^{-p-1}\text{log}^{q-\gamma }t+w\left( t\right) ,
\end{equation*}
where $\gamma =0$ if $p<0,$ $\gamma =1$ if $p\geq 0,$ $a=a(p,q)$ is a
constant and
\begin{equation*}
\left| \frac{d^{j}}{dt^{j}}w\left( t\right) \right| <C_{j}|\frac{d^{j}}{%
dt^{j}}\left[ t^{-p-1}\text{log}^{q-\gamma -1}t\right] ,\text{ \ \ \ \ }%
t\rightarrow \infty ,\text{ \ \ \ }j=0,1,2...\text{ \ }.
\end{equation*}
\end{lemma}

\bigskip

\bigskip

\end{document}